\documentclass[12pt]{amsart}
\usepackage{amsmath,amsthm,amsfonts}
\begin{document} 
\newcommand{\B}{{\mathbb B}}
\newcommand{\C}{{\mathbb C}}
\newcommand{\N}{{\mathbb N}}
\newcommand{\Z}{{\mathbb Z}}
\renewcommand{\P}{{\mathbb P}}
\newcommand{\R}{{\mathbb R}}
\newcommand{\rc}{\subset}
\newcommand{\dimc}{\mathop{dim}_{\C}}
\newcommand{\Lie}{\mathop{Lie}}
\newcommand{\Auto}{\mathop{{\rm Aut}_{\mathcal O}}}
\newtheorem*{definition}{Definition}
\newtheorem*{Conjecture}{Conjecture}
\newtheorem*{SpecAss}{Special Assumptions}
\newtheorem{example}{Example}
\newtheorem*{remark}{Remark}
\newtheorem*{remarks}{Remarks}
\newtheorem{lemma}{Lemma}
\newtheorem{proposition}{Proposition}
\newtheorem{theorem}{Theorem}
\title{%
Large Discrete Sets in Stein manifolds
}
\author {J\"org Winkelmann}
\begin{abstract}
We construct infinite discrete subsets of Stein manifolds with
remarkable properties. These generalize 
results of Rosay and Rudin on
discrete subsets of $\C^n$.
\end{abstract}
\address{%
J\"org Winkelmann \\
Universit\"at Basel\\
Mathematisches Institut \\
Rheinsprung 21 \\
CH-4051 Basel\\
Switzerland}
\email{jwinkel@member.ams.org\newline\indent{\itshape Webpage: }%
http://www.math.unibas.ch/\~{ }winkel/
}
\maketitle
\section{introduction}
Rosay and Rudin constructed examples of discrete subsets in $\C^n$
with rather remarkable properties (see \cite{RR}):
\begin{enumerate}
\item
There exists a discrete subset $D\subset\C^n$ such that the identity
map is the only automorphism of $\C^n$ stabilizing $D$.
\item 
There exists a discrete subset $D\subset \C^n$ such that for every
non-degenerate holomorphic map $F:\C^n\to\C^n$ the image $F(\C^n)$
has non-empty intersection with $D$.
\item
There are uncountably many inequivalent discrete subsets of $\C^n$
(where two discrete subsets $D_1$, $D_2$ are called equivalent if
there exists an automorphism $\phi$ of $\C^n$ such that
$\phi(D_1)=D_2$).
\end{enumerate}

Recently, Kaliman, using similar methods, showed the following
(\cite{Ka}).
\begin{enumerate}
\setcounter{enumi}{3}
\item
There exists a discrete subset $D\subset \C^n$ such that
$\C^n\setminus D$ is measure-hyperbolic.

\end{enumerate}
(Measure hyperbolicity is a notion introduced by Eisenman
in \cite{Eis}, see also \cite{KoB}.)

The goal of this article is to generalize these constructions to
complex manifolds other than $\C^n$. 
On a first glance it might appear that $\C^n$ is the most complicated
case and that therefore similar statements for arbitrary manifolds
should be easy corollaries of the statements for $\C^n$. However,
this is not the case. Indeed, there are simple
counterexamples to properties 1, 2 und 4.
\begin{proposition}[see ex.~\ref{ce-chr}, ex.~\ref{ce-nondeg},
ex.~\ref{ce-mh}]
\begin{enumerate}
\item
Let $V$ be a Stein manifold, $X=\P_n\setminus\{[1:0:0]\}$ with
$n=\dim V>0$ and
$D\subset X$ a discrete subset.

Then there exists a non-degenerate holomorphic mapping $F:V\to
X\setminus D$ and $X\setminus D$ is not measure hyperbolic.
\item
Let $D$ be a discrete subset of $X=\C\times\P_1$.

Then there exists a non-trivial automorphism $\phi$ of $X$ with
$\phi|_D=id_D$.
\end{enumerate}
\end{proposition}

Thus it is necessary to impose some conditions in order to generalize
these statements to complex manifolds other than $\C^n$.

Analyzing the methods of Rosay and Rudin, we came to the conclusion
that the crucial property of $\C^n$
is the existence of an exhaustion by
taut relatively compact domains, namely balls.
Given such an exhaustion it is possible to replace arguments using the
affine-linear structure of $\C^n$ by arguments based on the tautness
of these relatively compact domains, i.e.~normal family arguments.
Manifolds admitting such an exhaustion are called ``weakly Stein'' in
this article. Here is the precise definition:
\begin{definition}\label{def-ws}
A complex manifold $X$ is called \begin{em} weakly Stein \end{em} if
there exists a 
plurisubharmonic exhaustion function $\rho:X\to\R^+_0$
 such that 
\begin{enumerate}
\item
For every $c>0$ the set $B_c=\{x\in X:\rho(x)<c\}$ is taut.
\item
$S_c=\{x\in X:\rho(x)=c\}$ has empty interior for
all $c$.
\end{enumerate}
\end{definition}
In the sequel, we will simply say that $(Y,\rho)$ is weakly Stein in
order to indicate that $Y$ is a weakly Stein manifold and $\rho$ a
plurisubharmonic exhaustion function on $Y$ with the properties
required in the definition.
Every Stein manifold is weakly Stein (see \S\ref{sect-ws}).

We prove the generalization of the results of Rosay, Rudin and
Kaliman.

\begin{theorem}
Let $X$ be a weakly Stein manifold,
and let $V$ be an irreducible affine variety with $\dim V=\dim X>0$.

Then the following assertions hold:
\begin{enumerate}
\item
There exists a discrete subset $D\subset X$ such that $\phi(D)\ne D$
for every automorphism $\phi$ of $X$ except the identity.
\item
There exists a discrete subset $F\subset X$ such that $F(V)\cap
D\ne\{\}$ for every non-degenerate holomorphic map $F:V\to X$.
\item
There exist uncountably many inequivalent discrete subsets of $X$.
\item
There exists a discrete subset $D\subset X$ such that $X\setminus D$
is measure hyperbolic.
\end{enumerate}
\end{theorem}

The first assertion is theorem~\ref{th-main}, 
the second theorem~\ref{th-unavoid},
the third theorem~\ref{th-inequiv} 
and the last one is theorem~\ref{th-mh}.

\section{Weakly Stein manifolds}\label{sect-ws}
Let $X$ be a Stein manifold. Then there exists an 
embedding $i:X\hookrightarrow\C^n$ and
$\rho(x)\stackrel{def}{=}||i(x)||^2$ defines a plurisubharmonic
exhaustion function on $X$. For every $c>0$
the open subset $B_c=\{x\in X:\rho(x)<c\}$ is biholomorphic
to a closed submanifold of a ball in $\C^n$ and therefore
taut. It is obvious that $S_c=\{x\in X:\rho(x)= c\}$
has empty interior.
Thus a Stein manifold is necessarily weakly Stein.

More generally we have the following:
 
\begin{lemma}
Let $X$ be a complex manifold. Then the following four conditions
are equivalent:
\begin{enumerate}
\item
The complex manifold $X$ is both weakly Stein and holomorphically
convex.
\item
$X$ is holomorphically convex and every compact complex analytic
subspace of $X$ is hyperbolic (in the sense of Kobayashi).
\item
There exists a Stein space $Z$ and a proper surjective holomorphic map
$\pi:X\to Z$ such that all fibers are hyperbolic.
\item
For every infinite discrete subset $D\subset X$ and for every non-constant
holomorphic map
$g:\C\to X$ there exists a holomorphic function $f$ on $X$ such that
$f$ is unbounded on $D$ resp. $g(\C)$.  
\end{enumerate}
\end{lemma}
\begin{proof}
First recall that a complex manifold $X$ is holomorphically
convex if and only if for every discrete subset $D$ there exists a
holomorphic function $f$ such that $f|_D$ is unbounded and that
this property is equivalent to the existence of a proper, surjective,
connected holomorphic map from $X$ onto a Stein space $Z$.
Using this, it is easy to check
that each of the four properties implies that $X$ is
holomorphically convex and that there is a proper, surjective,
connected holomorphic map $\pi:X\to Z$ onto a Stein space $Z$.

Now let $W$ be a compact complex analytic subspace of a weakly
Stein manifold $(X,\rho)$.
Then $W$ is contained in some $B_c=\{x\in X:\rho(x)<c\}$,
because $W$ is compact. Thus $W$ is a closed subspace of a
taut manifold and therefore taut and in particular hyperbolic.
Hence compact complex analytic
subspaces of weakly Stein manifolds are hyperbolic
and therefore $(1)$ implies $(2)$.

Since $\pi:X\to Z$ is proper, its fibers are compact.
Thus $(2)$ implies $(3)$.

Assume $(3)$. Since $Z$ is Stein, it admits a strictly
plurisubharmonic exhaustion function $\tau$. For each $c>0$ the open
set $U_c=\{z\in Z:\tau(z)<c\}$ is complete hyperbolic.
Since $\pi$ is proper, the relative version of Brody's theorem
(see \cite{Lang})
implies that $B_c=\{x\in X:\pi(\tau(x))<c\}$ is complete hyperbolic
and therefore taut for every $c\in\R$.
Hence $(3)$ implies that $(X,\rho\circ\pi)$ is weakly Stein.
Furthermore, hyperbolicity of the $\pi$-fibers implies that
$\pi\circ g:\C\to Z$ can not be constant for any non-constant
holomorphic map $g:\C\to X$. Since $Z$ is Stein, it follows that
for every non-constant holomorphic map $g:\C\to X$ there exists
a non-constant holomorphic function $f_0$ on $Z$ such that 
$f_0\circ\pi\circ g$ is non-constant on $\C$.
Thus $(3)$ implies $(4)$.

Finally observe that
condition $(4)$  implies that there is no non-constant
holomorphic map from $\C$ into a fiber of $\pi$. By Brody's theorem 
(\cite{Br}) it follows that the fibers of $\pi$ are hyperbolic.
Hence $(4)$ implies $(3)$.
\end{proof}

There are also weakly Stein manifolds which are not
holomorphically convex.
Namely, let $M$ be a projective hyperbolic manifold and $L$ a flat line
bundle on $M$ such that none of the $L^n$ with $n\in\Z\setminus\{0\}$
is holomorphically trivial. Let $||\cdot||$ be a compatible hermitian
metric. Then $||\cdot||$ is a plurisubharmonic exhaustion function
having all the required properties.
Thus $L$ is weakly Stein.
On the other hand there are no non-constant holomorphic functions on
$L$. Hence the non-compact complex manifold $L$ is not  
 holomorphically convex.

\section{Preparations}
We begin by translating resp.~slightly generalizing two technical
results of \cite{RR}.
For this we utilize the following elementary lemma.
\begin{lemma}\label{le-2}
Let $(Y,\rho)$ be a weakly Stein  manifold with the corresponding
exhaustion function, $X$ a connected 
complex manifold of the same dimension and
$F:X\to Y$ be a non-degenerate (i.e.~having maximal rank at some
point) holomorphic map.
Assume that $F(X)\subset\{y\in Y:\rho(y)\le c\}$.

Then $F(X)\subset\{y\in Y:\rho(y)< c\}$.
\end{lemma}
\begin{proof}
By the maximum principle for plurisubharmonic functions it follows
that either the conclusion of the lemma holds or 
$F(Y)\subset\{y\in Y:\rho(y)=c\}=S_c$. However, the latter possibility is
excluded by the assumptions that $F$ is non-degenerate and that $S_c$
has empty interior.
\end{proof}
The following is a translation of lemma 4.3.~of \cite{RR}.
\begin{lemma}\label{le-3}
Let $U_2$ be a complex manifold, $U_1\subset U_2$ a relatively compact
connected open subset and
$C\subset U_1$ a compact subset. 
Furthermore let $(Y,\rho)$ be weakly Stein, $0<C_1<C_2$, $V_i=\{y\in Y:\rho(y)<C_i\}$
 and
$K\subset V_1$ be a
compact subset. Moreover fix hermitian metrics on $U_2$ and $V_2$
and a real constant $c>0$.

Let $\Gamma$ denote the set of all holomorphic mappings $F$ from $U_2$ to
$V_2$ such that $F(C)\cap K\ne\emptyset$ and $||DF||_C\ge c$
(where $||DF||_C$ is the supremum over all $||DF||_x$ for $x\in C$,
calculated with respect to the choosen Hermitian metrics).
Fix also a natural number $m\in\N$.

Then there exists a finite set 
$E=E(C,U_1,U_2;K,V_1,V_2;c,m)\subset\partial V_1$ with the
following property:

If $F\in\Gamma$ then either $F(U_1)\subset V_1$ or $F(U_2)\cap E$
contains at least $m$ points.
\end{lemma}
\begin{proof}
For the case $m=1$ the proof goes exactly as in \cite{RR}
(using lemma~\ref{le-2} instead of the maximum principle).
Observe that the method of construction allows one to
require that $E$ is contained in a fixed dense
subset of $\partial V$. Hence, after having constructed one such set $E$
we may construct another one, disjoint to the first one.
Continuing in this way, we may construct $m$ such sets and the union
of all these $m$ sets is the desired set $E$.
\end{proof}
We now slightly generalize lemma 4.4.~of
\cite{RR}.
\begin{lemma}\label{le-4}
Let $U_2$ be a complex manifold, $U_1\rc U_2$ a relatively compact
connected open subset, $C\subset U_1$ a compact subset, $c>0$,
$(Y,\rho)$ 
be a weakly Stein 
manifold, $r_n$ an strictly increasing sequence of positive numbers with 
$\lim r_n=+\infty$,
$V_n=\{y\in Y:\rho(y)<r_n\}$, 
$K\subset V_1$ a compact subset and assume that both
$U_2$ and $Y$ are equipped with hermitian metrics.
Let $\Gamma$ denote the set of all holomorphic mappings $F:U_2\to Y$
such that $F(C)\cap K\ne\emptyset$ and $||DF||_C\ge c$.
Let $m\in\N$.

Then there exists a discrete subset 
$E'=E'(C,U_1,U_2;K,V_1,Y;c,m)$ of $Y\setminus
V_1$ such that for every $F\in\Gamma$ either $F(U_1)\subset V_1$ or
$F(U_2)\cap E'$ contains at least $m$ points.
\end{lemma}
\begin{proof}
Choose an increasing sequence of connected
open subsets $W_i\subset U_2$ ($i\in\N)$) such
that
\begin{enumerate}
\item
Every $W_k$ is relatively compact in $W_{k+1}$
\item
$W_1=U_1$
\item
The union of all $W_k$ is relatively compact in $U_2$.
\end{enumerate}
Let
$
E'=\cup_{k\ge 1} E_k$ with 
\[
E_k=E'(C,W_k,W_{k+1};K,V_k,V_{k+1};c,m)
\]
Evidently $E'$ is discrete and $E'\cap V_1=\emptyset$.
Since $\cup_k W_k$ is relatively compact in $U_2$, for every
$F\in\Gamma$ there exists a number $n\in\N$ such that $F(W_k)\subset
V_n$ for all $k$. In particular $F(W_n)\subset V_n$.
Let us assume that $\#\left(F(U_2)\cap E'\right)<m$.
Clearly $\#\left(F(U_2)\cap E'\right)<m$ implies 
$\#\left(F(W_n)\cap E_{n-1}\right)<m$.
Hence lemma~\ref{le-3} yields $F(W_{n-1})\subset F(V_{n-1})$.
Going down by induction, we end up with the conclusion
$F(U_1)=F(W_1)\subset F(V_1)$ for all $F\in\Gamma$ with
$\#\left(F(U_2)\cap E'\right)<m$.
\end{proof}

\begin{proposition}\label{prop-1}
Let $X$ be a hermitian complex manifold with an increasing sequence of
open connected subsets $U_k$ such that each $U_k$ is relatively
compact in $U_{k+1}$.
Furthermore let $(Y,\rho)$ 
be a weakly Stein 
manifold equipped with a hermitian metric, 
$r_n$ an strictly increasing sequence of positive numbers with $\lim r_n=+\infty$,
$V_n=\{y\in Y:\rho(y)<r_n\}$ 
and $\dimc(X)=\dimc(Y)$.

Then there exists a discrete subset $D\subset Y$ and a mapping
$\alpha:\N\times\N\times\R^+\to\N$ with the following property:

Let $F:X\to Y$ be a non-degenerate holomorphic map such that
there is a number $M\in\N$ such that
$\#\left( F(U_{n+1})\cap D\right)\le M+\#\left( V_n\cap D\right)$ 
for all
$n\in \N$.

Then
$F(U_n)\subset V_n$ for all
$n\ge\alpha(M,q,c)$
where $q=\min\{k\ge 2:F(\overline{U_k})
\cap \overline{V_k}\ne\emptyset\}$ and 
$c=||DF||_{\overline{U_q}}$.
\end{proposition}
\begin{proof}
Choose discrete subsets $D_k\subset Y\setminus V_k$
by 
\[
D_k= E'(\overline{U_{k-1}},U_k,U_{k+1};\overline{V_{k-1}},V_k,Y;
1/k,k+1+\#\left((\cup_{i<k}D_i)\cap V_k\right))
\]
and $D=\cup_k D_k$. Note that $D$ is discrete, because $D_k\subset
Y\setminus V_k$.

Let $F:X\to Y$ be a holomorphic map which is non-degenerate, i.e.~%
$||DF||\not\equiv 0$. 
Let $q\in\N$ be choosen in such a way that
$F(\overline{U_{q}})\cap\overline{V_q}\ne\emptyset$. 
Furthermore $c=||DF||_{\overline{U_{q}}}>0$. 
Finally assume that there is a number $M\in\N$ such that
$\#\left( F(U_{n+1})\cap D\right)\le M+\#\left( V_n\cap D\right)$ holds
for all $n\in\N$.
Then for all $l>M$ we obtain
\[
\begin{split}
&\#\left(F(U_{l+1})\cap D_l\right)\\
\le\ 
&\#\left(F(U_{l+1})\cap D\right)\\
\le\ 
&M+1+\#\left( V_l\cap D\right)\\
<\ 
&l+1 + \#\left( V_l\cap D\right)\\
\end{split}
\]
From the construction of $D$ it then follows that
$F(U_l)\subset V_l$ for all $l>\max(M,q+1,1/c)$.
Hence the statement of the proposition is fulfilled if we choose a
map
$\alpha$ in such a way that $\alpha(M,q,c)\ge\max(M,q+1,\frac{1}{c})$
for all $(M,q,c)\in\N\times\N\times\R^+$.
\end{proof}
\section{Rigid sets}
For this section we fix the following presumptions.
\begin{SpecAss}
\begin{enumerate}
\item
$(X,\rho)$ is weakly Stein (as defined in Def.~\ref{def-ws}).
\item
$(r_n)_{n\in\N}$ is an strictly increasing sequence of positive numbers with 
$\lim r_n=+\infty$.
\item
$U_n=\{x\in X:\rho(x)<r_n\}$.
\item
$D$ is a discrete subset of $X$.
\item
$\alpha:\N\times\N\times\R^+\to\N$ is a map such that for every
$f\in Aut(X)$ the condition
\[
\# f(U_n\cap D)\le m+\#(U_n\cap D)\ \forall n\in\N
\]
together with $f(\overline{U_q})\cap\overline{U_q}\ne\emptyset$
implies that
\[
f(U_n)\subset f(U_{n+1})\ \forall
n\ge\alpha(m,q,||DF||_{\overline{U_q}}).
\]
\end{enumerate}
\end{SpecAss}
Prop.~\ref{prop-1} above implies 
that for every weakly Stein manifold $(X,\rho)$
and every such sequence $(r_n)$ there exists a discrete subset $D$
and a map $\alpha$ such that these special assumptions are fulfilled.

\begin{proposition}
Under the special assumptions let $D'$ be a  subset of $X$
such that $(D\setminus D')\cup (D'\setminus D)$ is finite.

Then $Aut(X,D')=\{\phi\in Aut(X):\phi(D')=D'\}$ is a
real Lie group (with countably many connected components).
\end{proposition}
\begin{proof}
The automorphism group $Aut(X)$ embeds into $Hol(X,X)\times Hol(X,X)$ 
via $\phi\mapsto(\phi,\phi^{-1})$. We equip $Aut(X)$ with the relative
topology corresponding to this embedding ($Hol(X,X)$ being topologized
with the compact-open-topology).
The topology on $Aut(X,D')=\{\phi\in Aut(X):\phi(D')=D'\}$ is induced
by the inclusion $Aut(X,D')\subset Aut(X)$.

Consider
\[
W=\{\phi\in Aut(X,D'):
\phi(\overline{U_1})\subset\overline{U_2}\text{ and }
\phi^{-1}(\overline{U_1})\subset\overline{U_2}
\}
\]
We will prove that $W$ is compact.
For this note that
\[
||D\phi||_p = ||D\phi^{-1}||_{\phi(p)}^{-1}
\]
for all $p\in X$, $\phi\in Aut(X)$. 
Now $\phi(\overline{U_2})\cap\overline{U_2}\ne\emptyset$ for $\phi\in
W$.
Therefore
\[
||D\phi||_{{\overline{U}}_2} \cdot
||D(\phi^{-1})||_{{\overline{U}}_2} \ge 1
\]
for all $\phi\in W$.
Let $(\phi_l)_l$ be a sequence in $W$.
By passing to a subsequence and possibly replacing $(\phi_l)_l$ by the
sequence of the respective inverse maps $(\phi_l^{-1})_l$, we may assume
that $||(D\phi_l)||_{{\overline{U}}_2}\ge 1$ for all $l$.
Furthermore observe that
\[
\begin{split}
&\#\phi_l(U_n\cap D) \\
\le\ & \#(D\setminus D') + \#\phi_l(U_n\cap D')\\
=\ & \#(D\setminus D') + \#(U_n\cap D')\\
\le\ & \#(D\setminus D') + \#(D'\setminus D) + \#(U_n\cap D).\\
\end{split}
\] 
Hence
$\phi_l(U_n)\subset U_{n+1}$ for all $l$ and $n$ with 
$n\ge N=\alpha(\#(D\setminus D')+\#(D'\setminus D),2,1)$ by special
assumption $(5)$. Note that $N$ does not depend on $l$.
Therefore a subsequence of $(\phi_l)_l$ converges to a holomorphic map
$F:X\to X$.
This convergence implies
$\limsup_l||D\phi_l||_{{\overline{U}}_2}=\zeta<\infty$.
Hence 
$\phi_l^{-1}(U_n)\subset U_{n+1}$ for
$n\ge\alpha(\#(D\setminus D')+\#(D'\setminus D),2,\zeta)$
and we obtain that the sequence of inverse maps $\phi_l^{-1}$
converges as well. But now $\lim\phi_l=F$ and $\lim\phi_l^{-1}=G$
imply that $F\circ G=id_X$. Hence $F,G\in Aut(X)$.
Furthermore it is clear that $F(D')\subset D'$ and $G(D')\subset D'$.
It follows that $F,G\in Aut(X,D')$.
Thus $W$ is compact. Since $W$ contains an open neighbourhood of
$id_X$ in $Aut(X,D')$, it follows that $Aut(X,D')$ is a locally
compact topological group.
By a result of Bochner and
Montgomery
\cite{BM} it follows that $Aut(X,D')$ is a real Lie group.

Finally note that $X$ having countable base of topology implies that
$Hol(X,X)$ and consequently $Aut(X,D')$ have likewise topologies with
countable base. Since $Aut(X,D')$ is a Lie group, its connected
components are open. Hence the countability of the base of the topology
implies that $Aut(X,D')$ has only countably many connected components.
\end{proof}
\begin{remark}
This evidently implies that the isotropy groups
$Aut(X,D')_p\allowbreak
=\{g\in Aut(X,D'):g(p)=p\}$ are likewise Lie groups with
countably many connected components. By arguments similar to the ones
used in the above proof one may  show in addition that $Aut(X,D')_p$ is
actually compact for all $p\in X$.
\end{remark}

\begin{lemma}
Let $K$ be a  real Lie group 
(with countably many connected components)
acting effectively on a connected complex
manifold by biholomorphic transformations.

Then there exists a finite set $S$ such that
$k\in K$,
 $k(x)=x$ for all $x\in S$ implies $k=e$.
\end{lemma}
\begin{proof}
Let $n=\dim K$ and choose $x_1,\ldots x_n$ in such a way that
\[
\dim\{g\in K:g(x_i)=x_i\ \forall 1\le i \le k\}\le n-k.
\]
Then $\Gamma=\{g\in K:g(x_i)=x_i\ \forall 1\le i\le n\}$ 
is discrete. For every $\gamma\in\Gamma$ the set $X^\gamma=
\{x\in X:\gamma(x)=x\}$ is a closed analytic subset and therefore of
measure zero.
Thus $\cup_{\gamma\in\Gamma} X^\gamma$ is of measure zero, too.
Now choose $S=\{x_1,\ldots,x_n,y\}$ with $y\not\in\cup_\gamma
X^\gamma$.
\end{proof}
\begin{proposition}
Under the special assumptions there exists a finite set 
 $S\subset X$ such that 
the following properties are fulfilled:
\begin{enumerate}
\item
For every discrete subset $D'$ of $X$ with $D\cup S\subset D'$
the group $Aut(X,D')$ is countable.
\item
If $\phi\in Aut(X)$ and $\phi(x)=x$ for all $x\in D\cup S$, then
$\phi=id_X$.
\end{enumerate}
\end{proposition}

\begin{proof}
Choose $p\in D$ and let $K=Aut(X,D)_p$, i.e.~%
$K=\{\phi\in Aut(X):\phi(D)=D\text{ and }\phi(p)=p\}$.
Then $K$ is a  real Lie group acting effectively on the
manifold
$X\setminus D$.
By the preceding lemma there is a finite set $S\subset X\setminus D$
such that $k|_S=id_S$ implies $k=id_X$ for $k\in K$.

Now let $D'$ be a discrete subset of $X$ with $D\cup S\subset D'$
and let $G$ be the connected component of $id_X$ in
$Aut(X,D')$.
Since $G$ is connected and stabilizes $D'$, it
must act trivially on $D'$. Hence $G\subset K$ and therefore
$G=\{id\}$ by the choice of $S$.
It follows that $Aut(X,D')$ is countable.

Finally, if $\phi\in Aut(X)$ with $\phi|_{D\cup S}=id_{D\cup S}$,
then $\phi\in K$ and thus $\phi=id_X$ by the choice of $S$.
\end{proof}

\begin{theorem}\label{th-main}
Let $X$ be a weakly Stein manifold.

Then there exists a discrete subset $\Lambda$ with the property that
$\phi\in Aut(X)$ with $\phi(\Lambda)=\Lambda$ implies $\phi=id_X$.
\end{theorem}

\begin{proof}
Choose a finite set $S$ as in the above proposition. Recall that the
choice of $S$ is generic, i.e.~there is no problem in finding another
finite subset $S'$ with $S\cap S'=\emptyset$ enjoying the same
properties as $S$.
 
Let $D_1=D\cup S\cup S'$.
Then for every $x\in D$ we have the conclusion that
$Aut(X,D_1\setminus\{x\})$ is countable and
that $\phi|_{D_1\setminus\{x\}}=id$
implies $\phi=id_X$.

Let $\Omega$ be the set of all $\omega\in X\setminus D_1$ 
for which there exists
an automorphism $\phi\in Aut(X)$ such that $\phi$
stabilizes $D_1\cup\{\omega\}$ but $\phi(\omega)\ne\omega$.
For every $\omega\in\Omega$ let us fix one such $\phi$, henceforth
called $\phi_\omega$.
For every $(a,b)\in D_1\times D_1$ let $\Omega(a,b)=
\{\omega\in\Omega:\phi_\omega(\omega)=a\text{ and }
\phi_\omega(b)=\omega \}$
and for every non-empty 
$\Omega(a,b)$ fix an element $\eta(a,b)\in\Omega(a,b)$.
Then we obtain a map from $\Omega(a,b)$ to $Aut(X,D_1\setminus\{b\})$ by
\[
\omega \mapsto \phi_\omega^{-1}\circ\phi_{\eta(a,b)}
\]
This map is injective, because 
\[
(\phi_\omega^{-1}\circ\phi_{\eta(a,b)})|_{D_1\setminus\{b\}}
=(\phi_{\tilde\omega}^{-1}\circ\phi_{\eta(a,b)})|_{D_1\setminus\{b\}}
\]
implies $\phi_\omega=\phi_{\tilde\omega}$ which in
turn implies 
\[
\omega=\phi_\omega(b)=\phi_{\tilde\omega}(b)=\tilde\omega
\]

It follows that $\Omega(a,b)$ is countable for all $(a,b)\in D_1\times
D_1$.
Therefore $\Omega$, being the union of all $\Omega(a,b)$, is countable
as well.

For every $\phi\in Aut(X)$ the fixed point set $X^\phi$ is an analytic
subset of $X$.
Since $Aut(X,D_1)$ is countable, it follows that 
\[
\Omega'=\{x\in X:\exists\phi\in Aut(X,D_1)\setminus\{id\}:\phi(x)=x\}
\]
is a set of measure zero.

Now $D_1$, $\Omega$ and $\Omega'$ are sets of measure zero.
Hence there
 is a point $p\in X$ contained in none of these three sets.
Let $\Lambda=D_1\cup\{p\}$.
We claim that $Aut(X,\Lambda)=\{id\}$.
Indeed, $p\not\in\Omega$ implies that there is no $\phi\in Aut(X,\Lambda)$
 with $\phi(p)\ne p$. Hence $Aut(X,\Lambda)\subset Aut(X,D_1)$.
But $p\not\in\Omega'$ implies that $\phi(p)\ne p$ for all $\phi$ in
 $Aut(X,D_1)$ except $id_X$. It follows that $Aut(X,\Lambda)=\{id\}$.
\end{proof}

The assumption of $X$ being weakly Stein was used to derive that
certain subgroups of $Aut(X)$ are small in a certain sense.
Therefore we may replace the assumption of $X$ being weakly Stein
by the assumption that $Aut(X)$ is sufficiently small, e.g. a Lie
group.
One can check that the arguments used in the proof of the theorem
above also imply the result stated below.
\begin{proposition}
Let $X$ be a complex manifold, $\dim(X)>0$,
 such that $Aut(X)$ is a real Lie group.

Then there exists a discrete subset $D\subset X$ such that 
$\phi(D)\subset D$ implies $\phi=id_X$ for $\phi\in Aut(X)$.
\end{proposition}
For certain complex manifolds $X$ it is true that $Aut(X)$ must be a
real Lie group, e.g.~this is true for compact complex manifolds
as well as for bounded domains in $\C^n$ and other hyperbolic
complex manifolds.

However, in general the automorphism group of a complex manifold
does not need to be a Lie group. E.g., if $X$ is a connected
complex manifold
 and $Y$ is a complex
manifold with a complex one-parameter-group of automorphisms
$\psi_t\in Aut(Y)$ then for a every holomorphic function $f$ on $X$ a
holomorphic automorphism of $X\times Y$ is given by
\[
(x,y)\mapsto (x,\psi_{f(x)}(y)).
\]
This implies that $Aut(X,Y)$ can not be a Lie group unless
every holomorphic function on $X$ is constant.

In this spirit we derive the following example.
\begin{example}\label{ce-chr}
Let $Z$ be a compact complex manifold on which a complex Lie group G
acts effectively, $Y$ a Stein manifold and $X=Y\times Z$
(with $\dim(G),\dim(Y)>0$). Then for every discrete
subset $E\subset X$ there exists a non-trivial automorphism $\phi\in Aut(X)$
with $\phi|_E=id_E$.
\end{example}
\begin{proof}
Let $\pi_2:X\to Y$ denote the projection onto the second factor.
This map is proper, hence $\pi_2(E)$ is discrete in $Y$. Now simply
choose a non-zero
holomorphic function $f\in {\mathcal O}(Y)$ vanishing on
$\pi_2(E)$ and an element $v\in\Lie(G)$ and let $\phi$
be the automorphism of $X$ defined by $(z,y)\mapsto
(\exp(vf(y))(z),y)$.
\end{proof}
\section{Inequivalent discrete subsets}
The methods used in the preceding sections also may be employed to
show the following.
\begin{theorem}\label{th-inequiv}
Let $X$ be a weakly Stein manifold.

Then there exist uncountably many pairwise inequivalent infinite discrete
subsets of $X$.
\end{theorem}
Here two discrete subsets $D$, $D'$ of $X$ are called equivalent if
there exists an automorphism $\phi\in Aut(X)$ such that $\phi(D)=D'$.
\begin{proof}
Let $D_1$, $\Omega$ and $\Omega'$ be as in the proof of theorem~\ref{th-main}.
Let $\Sigma=D_1\cup\Omega\cup\Omega'$.
For $x,y\in X\setminus\Sigma$ we define an equivalence relation by setting
$x\sim y$ iff there exists an automorphism $\phi\in Aut(X)$ such that
$\phi(D_1\cup\{x\})=D_1\cup\{y\}$. 
For $x\in X\setminus\Sigma$ and $a\in D_1$ let
\[
E_{a}(x)= \{ y\in X\setminus\Sigma:
\exists \phi\in Aut(X):\phi(D_1\cup\{x\})=D_1\cup\{y\},
\phi(a)=y \}
\]
and 
\[
E^*(x)=\{y\in X\setminus\Sigma:\exists\phi\in Aut(X,D_1):\phi(x)=y\}
\]
Recall that $Aut(X,D_1)$ is countable. 
Hence $E^*(x)$ is also countable.
Let $y,\tilde y\in E_{a}(x)$ and let $\phi$, $\tilde\phi$ be
corresponding
elements in $Aut(X)$. Then $\psi=\tilde\phi\circ(\phi^{-1})\in Aut(X)$ with
$\psi(y)=\tilde y$ and $\psi(D_1)=D_1$. Thus countability of
$Aut(X,D_1)$ also implies that all the sets $E_a(x)$ are countable.
Therefore  the sets $E_a(x)$ and $E^*(x)$ are all countable as well.
It follows that the
equivalence classes of $\sim$ in $X\setminus\Sigma$ are countable.
Since $X\setminus\Sigma$ is not countable, this implies
 that there are
uncountably many equivalence classes. 
\end{proof}
In fact this result can be proved for many other spaces as well.
For instance there is the following proposition.
\begin{proposition}
Let $X$ be a complex manifold. Assume that there exists a
non-constant holomorphic function $f$ on $X$.

Then there exist two inequivalent infinite discrete
subsets of $X$.
\end{proposition}
\begin{proof}
Let $r=\sup\{|f(x)|:x\in X\}$. If $r<\infty$, then choose a sequence
$x_n$
in $X$ with $\lim |f(x_n)|=r$ and let $a$ be an accumulation point
 of $f(x_n)$. In this case $1/(f-a)$ is an unbounded holomorphic
function on $X$. Thus it is clear that for some sequences $x_n$
in $X$ there is a holomorphic function $f$ such that $|f(x_n)|$ is 
unbounded. If there are sequences on which every holomorphic function
is bounded, we found two inequivalent infinite discrete subsets.
Thus we may assume that there is no such sequence, i.e., we may
assume that $X$ is holomorphically convex and we may consider
the \begin{em} Remmert-reduction $\pi:X\to Z$ \end{em}
(where $\pi$ is connected, surjective and proper and $Z$ is Stein).

If $\pi$ is not injective, there is an infinite discrete subset
$S\subset X$ such that $\pi|_S$ is not injective and such an infinite
discrete subsets is not equivalent to any infinite discrete subset
$S'\subset X$ for which $\pi|_{S'}$ is injective.
As a consequence, we may assume that $\pi$ is injective, i.e.,
that $X\simeq Z$, i.e., that $X$ is Stein.

Now the assertion follows from the preceding proposition.
\end{proof}

\begin{Conjecture}
For every  complex space $X$ with $\dim(X)>0$
there exists a uncountable family of
inequivalent discrete subsets
of the same cardinality.
\end{Conjecture}

For non-compact spaces this conjecture is motivated by the above
considerations. For a compact complex space $X$, the automorphism
group $\Auto(X)$ is a finite-dimensional complex Lie group $G$.
Then, for every natural number $n$ with $n\dimc(X)>\dimc(G)$ the induced
action of $\Auto(X)$ on the $n$-fold symmetric product $S^n(X)$ can
not have an open orbit. It follows that for such a number $n$ there
are uncountably many inequivalent subsets of $X$ with cardinality $n$.
Hence the conjecture holds for compact complex spaces.

\section{Unavoidable sets}
Let $X$, $Y$ be complex manifolds. A subset $S\subset Y$ is called
unavoidable for a family $\mathcal H$ of holomorphic mappings from $X$ to
$Y$ if $f(X)\cap S\ne\emptyset$ for every $f\in\mathcal H$.
Given a pair of complex manifolds $(X,Y)$ of the same dimension 
let $Hol^*(X,Y)$ denote the set of non-degenerate holomorphic maps.
Clearly, $Hol^*(X,Y)$-unavoidable discrete subsets do not exist for 
arbitrary pairs
of complex manifolds $(X,Y)$.
For instance, assume that $X$ is a  bounded domain in a Stein
manifold. Then there exists a non-degenerate holomorphic map from $X$
to the unit ball $B_n$ in $\C^n$. As a consequence there is a
non-degenerate holomorphic map from $X$ to $Y$ with $F(X)\cap
D=\emptyset$ for every complex manifold $Y$ and every 
closed analytic subset
$D\subset Y$.

It is also necessary to put some conditions on the target manifold,
even if it is non-compact.
\begin{example}\label{ce-nondeg}
Let $X$ be a Stein manifold of dimension $n$, $p\in\P_n$,
$Y=\P_n\setminus\{p\}$ and $D\subset Y$ a discrete subset.

Then there exists a non-degenerate holomorphic map $F:X\to Y$ with
$F(X)\cap D=\emptyset$.
\end{example}
\begin{proof}
Let $H\subset\P_n$ be a hyperplane which does not intersect
$D\cup\{p\}$. After changing coordinates, we may assume
that
$p=[1:0:\ldots:0]$ and $H=\{[x_0:\ldots:x_n]:x_0=0\}$.
Let $\phi\in Aut(\P_n)$ defined by 
\[
\phi:[x_0:\ldots:x_n]\mapsto
[x_0:\frac{1}{2}x_1:\ldots:\frac{1}{2}x_n].
\]
Now for any open neighbourhood $W$ of $p$ in $\P_n$ the complement
$\P_n\setminus W$ is compact and $D\cap(\P_n\setminus W)$ is finite.
Hence for any open neighbourhood $W$ of $p$ in $\P_n$
there exists a number $N(W)$ 
such that $\phi^{N(W)}(D)\subset W$.
It follows that there exists a biholomorphic map $\zeta$ from $\C^n$
onto some open neighbourhood of $p$ in $\P_n$ such that $\zeta(\B_n)$
contains all of $D$.
Next choose a non-degenerate holomorphic map $H:X\to \C^n$ and a
non-degenerate holomorphic map $G:\C^n\to\C^n$ such that $G(\C^n)$ does
not intersect the unit ball in $\C^n$ (such a map exists by the
Fatou-Bieberbach construction).
Then $F=\zeta\circ G\circ H:X\to \P_n\setminus\{p\}=Y$ has the desired
property $F(X)\cap D=\emptyset$.
\end{proof}

\begin{theorem}\label{th-unavoid}
Let $X$ be an irreducible, reduced affine variety and $(Y,\rho)$ be weakly
Stein.

Then there exists a discrete subset $D\subset Y$ such that $F(X)\cap
D\ne\emptyset$ for every non-degenerate holomorphic map 
$F:X\to Y$.
\end{theorem}
\begin{proof}
By the ``Noether Normalization Lemma'' there exists a finite ramified
covering $\pi:X\to\C^n$. Let $V_n=\{y\in Y:\rho(y)<n\}$ and
$U_n=\{x\in X:||\pi(x)||<e^{n^2}\}$.
By prop.~\ref{prop-1} there exists a discrete subset $D\subset Y$ such that
for every non-degenerate holomorphic map $F:X\to Y$ with $F(X)\cap
D=\emptyset$ there is a number $N$ such that $F(U_n)\subset V_n$ for
$n\ge N$. We claim that there can not exist such a map $F$.
Consider the plurisubharmonic function $F^*\rho$ on $X$ and define
a plurisubharmonic function $\xi$ on $\C^n$ by
\[
\xi(v)=\sum_{x\in\pi^{-1}(v)}F^*\rho(x)
\]
(This formula is to be understood literally where $\pi$ is
unramified, in the ramification locus multiplicities have to be taken
into account.)
Now $F(U_n)\subset V_n$ for $n>\!>0$ implies that
\[
\limsup_{||v||\to\infty} \frac {\xi(v)}{\log ||v||} = 0
\]
But $M(r)=\sup_{||v||\le r}\xi(v)$ is an increasing convex function
in $\log r$ (see e.g.~\cite{LG}, prop.~1.4.b).
For a real convex function $t\mapsto M(t)$ we have
\[
M(t)\le \frac{n-1}{n} M(0) + \frac{1}{n}M(nt)
\]
for all $t,n\in\R^+$.
Thus $\lim_{n\to\infty} \frac{M(nt)}{nt}=0$
implies $M(t)\le M(0)$ for all $t>0$.
This forces $\xi$ to be
constant.
Let $\xi\equiv C$. Recall that 
 $F^*\rho\ge 0$. It follows that $C\ge F^*\rho$. But plurisubharmonic
functions bounded from above on affine varieties are constant. Hence
$F^*\rho$ must be constant, which contradicts the assumption that
$F$ is non-degenerate.
\end{proof}

\section{Discrete sets with measure hyperbolic complements}

\begin{theorem}\label{th-mh}
Let $(Y,\rho)$ be a weakly Stein complex manifold.

Then there exists a discrete subset $D\subset Y$ such that
$Y\setminus D$ is measure--hyperbolic, i.e., the Kobayashi-Eisenman
pseudovolume on $Y\setminus D$ is a volume.
\end{theorem}

\begin{proof}
We will apply lemma~\ref{le-4}.
For this purpose, let $n=\dim Y$, $U_2=\B_n=
\{v\in\C^n:||v||<1\}$, $U_1=\frac{1}{2}\B_n$,
$C=\{(0,\ldots,0)\}\in\B_n$, $c=1$ and $Y_k=\{y\in Y:\rho(y)<k\}$.
Furthermore we endow $Y$ with an hermitian metric $h$ and denote the
corresponding volume form by $\omega$. 
Let the manifold $U_2=\B_n$ be endowed with
the euclidean metric.

For each $k\in\N$, we choose a discrete subset 
$E_k\subset Y\setminus Y_k$ by invoking lemma~\ref{le-4}
with $E_k=E'(C,U_1,U_2;\overline{Y_{k-1}},
Y_k,Y;1,1)$ and we define 
\[
D=\cup_{k\in\N} E_k
\]

Fix a point $p\in Y\setminus D$. Evidently $p\in\overline{ Y_{k-1}}$ 
for some $k\in\N$.
Let $f:\B_n\to Y\setminus D$ be a holomorphic map with $f(0)=p$.
Then $E_k\subset D$ implies $||Df||_0<1$ or
$f(\frac{1}{2}\B_n)\subset Y_k$.
For the  Kobayashi-Eisenman pseudovolume $\mu_{Y\setminus D}$
it follows that 
\[
\mu_{Y\setminus D}
\ge
\min\{ \omega, \frac{1}{2^n}\mu_{Y_k}
\}
\]
where $\mu_{Y_k}$ denote the Kobayashi-Eisenman pseudovolume
of the complex manifold $Y_k$.
Since $Y_k$ is taut and therefore measure-hyperbolic, it follows
that $\mu_{Y\setminus D}(p)>0$ for all $p\in Y$,
i.e., $Y\setminus D$ is measure-hyperbolic.
\end{proof}
\begin{remark}
A holomorphic map from $\C^n$ to a $n$-dimensional
measure hyperbolic complex manifold must be degenerate.
Thus we reproved theorem~\ref{th-unavoid} for the special case $X=\C^n$.
\end{remark}
\begin{example}\label{ce-mh}
Let $Y=\P_2\setminus\{[1:0:0]\}$. Then
$Y\setminus D$ is never measure hyperbolic for any discrete
subset $D\subset Y$.
To verify this, note that there exists a non-degenarate holomorphic
map from $\C^n$ to $Y\setminus D$ (see ex.~\ref{ce-nondeg}) and recall
that a holomorphic map from $\C^n$ to a measure hyperbolic manifold
is necessarily degenerate.
\end{example}

\end{document}